\documentclass[english,nofirstpagebreak,cropmarks]{studia-Hermann} 

\usepackage[T1]{fontenc} 
\usepackage[english]{babel}
\usepackage{cite}

\newtheorem{theorem}{Theorem}[section]
\newtheorem{defi}{Definition}[section]
\newtheorem{lemme}{Lemma}[section]

\newtheorem{cor}{Corollary}[section]
\newtheorem{proof}{Proof}[section]

\newtheorem{rem}{Remark}[section]

\journal{Studia Informatica Universalis.}{-}{-}

\title[\textit{The Volterra Models}]{Volterra Equations of the First kind with Discontinuous Kernels in the\\ Theory of Evolving Systems Control}
\author{D. N. Sidorov\fup{*}}
\address{\fup{*}Melentiev Energy Systems Institute\\
Russian Academy of Sciences\\
 130 Lermontov Street, 664033 Irkutsk\\
Irkutsk State University\\
20 Boulevard Gagarin,  664003 Irkutsk\\ 
dsidorov@isem.sei.irk.ru \\ [3pt]
}
\resume{Tapez votre r�sum� ici.}
\abstract{
The Volterra integral equations of the first kind
with piecewise smooth kernel are considered. Such equations appear in the theory of optimal control of the evolving systems. The existence theorems are proved. 
The method for constructing approximations of parametric families of solutions of such equations is suggested. The parametric family of solutions is constructed in terms of a logarithmic power
asymptotics.}
\motscles{\'Equations  int\'egrale de Volterra, noyau discontinu, approximations successives, contr\^
ole optimal, }
\keywords{Volterra  equations, discontinuous kernel, succesive approximations, optimal lifetime,
evolving systems}

\begin{document}
\maketitlepage


\section*{Introduction \label{sec1}}

The theory of integral models of evolving systems was first initiated  by V. Glushkov 
in the 70s of XX century. Readers may refer to the bibliorgaphy in \cite{1,2,3,4,5}. Such theory employs the Volterra integral equations of the first kind where
limit of integration is time function. The theory and numerical methods of such nonclassic equations were studied
in the monograph \cite{2} and applied for estimation and optimization of lifetime of electric
power systems (EPS) components \cite{4}.

In this paper we address  the following integral equation:
\begin{equation}
\int\limits_{0}^{t} K(t,s) x(s) ds = f(t), \,\,\,  0 < t \leq T,
\label{eq1}
\end{equation}
where kernel is defined as follows:
 \begin{equation}
    K(t,s) = \left\{ \begin{array}{ll}
         \mbox{$K_1(t,s), \,\,\,\,\,\,\,\,\,\,\,\, {0} \leq s < \alpha_1(t)$}, \\
         \mbox{$K_2(t,s), \,\, \alpha_1(t) \leq s < \alpha_2(t)$}, \\
         \mbox{\,\, \dots \,\,\,\,\,\,\,\,\,\,\,\,\, \dots \dots \dots } \\
         \mbox{$K_n(t,s), \,\, \alpha_{n-1}(t) \leq s \leq { t}$}, \\
        \end{array} \right. \label{eq2} 
         \end{equation}
         $$ 0 < \alpha_1(t) < \alpha_2(t) < \dots < \alpha_{n-1}(t) <   t, \,\, |\alpha_i^{\prime}(t)|<1. $$
Functions $K_i(t,s), i=\overline{1,n},  f(t)$ are continuous and sufficiently smooth, $f(0)=0.$

The Glushkov integral model of evolving systems \cite{1,4,5} is the special case of the Volterra integral equation (\ref{eq1}) -- (\ref{eq2}) when all the kernels
except $K_1(t,s)$ are zeros. 

First results in studies of the Volterra equations with discontinuous kernels were
formulated by G.C. Evans\cite{6} in the beginning of XX century.
Results in the spectral theory of integral operators with discontinuous kernels 
were obtained by A.P. Khromov in \cite{7}.
Asymptotic approximations of solutions of the Volterra equations 
of the first kind with analytical kernel $K(t,s)$ were constracted by N.A. Magnitsky \cite{8}. 

It is to be noted that solutions of the equations 
(\ref{eq1}) can have an arbitrary constants and can be unbounded for  $t \rightarrow 0.$
For example, if 
\begin{equation}
    K(t,s) = \left\{ \begin{array}{ll}
         \mbox{$1, \, 0\leq s < t/2 $}, \\
         \mbox{$-1, \, t/2 \leq s \leq t $}, \\
        \end{array} \right.
        \end{equation}
        $f(t)=t,$
        then equation (\ref{eq1})  has the solution
         $x(t) = c-\frac{\ln t}{\ln 2},$ where $c$ is constant.

        In this paper we employ results of the papers \cite{8, 9, 10, 11}
         in order to formulate the algorithm for construction of the continuous 
solutions of equation  (\ref{eq1}) for  
        $0 < t \leq T$ in the following form:
        \begin{equation}
        x(t) = \sum\limits_{i=0}^N x_i(\ln t) t^i + t^N u(t).
        \label{eq3}
        \end{equation}
Coefficients $x_i(\ln t)$
are constructed as polynomials on powers of $\ln t$ and they may depend on certain number of 
arbitrary constants.  $N$ defines the necessary smoothness of the functions $K_i(t,s), \, f(t).$ In this paper we propose an algorithm for construction of the function $u(t)$
in representation of the desired solution (\ref{eq3}) based on successive approximations
method which is uniformly converge on $[0,T].$  It is to be noted,
that logarithmic-power asymptotics have been efficiently employed for
solution of integral and differential equations in irregular cases\cite{10} -- \cite{12}, 
\cite{8}.

The paper is organized as follows. In Section 1 after the problem statement  we 
introduce the structure of  solution and  prove the existence theorem.
The method for the asymptotic approximations construction is suggested in
 Section 2. The main theorem is formulated and proved also in  Section 2.
Finally, concluding remarks are given.
 
 \section{The structure of solutions and  existence theorem}
 For sake of clarity let us suppose that  $\alpha_i(t) = \alpha_i t,$ where
 $0 = \alpha_0 < \alpha_1 < \cdots < \alpha_n = 1.$
 Let us introduce the condition\\
 {\bf A.} $K_n(t,t) \neq 0$ for $t \in [0,T]$
 and $N$ is selected to fulfill the following inequality
 \begin{equation}
 \max\limits_{0\leq t \leq T} |K_n(t,t)|^{-1} \sum\limits_{i=1}^n \biggl ( \alpha_i^{1+N}|K_i(t, \alpha_i t)| + \,\,\,\,\,\,\,\,\,\,\,\,\,\,\,\,\,\,\,\,\,\,\,\,\,\,\,\,\,\,\,\,\,\,\,\,\,\,\,\,\,\,\,\,\,\,\,\,
 \label{eq4}
 \end{equation}
\hfill $+\alpha_{i-1}^{1+N} |K_i(t, \alpha_{i-1} t)|   \biggr)
 \leq 1+q, $\\
where $q<1,$ $\alpha_0=0, \alpha_n=1.$
Condition {\bf A} is fulfilled for large enough $N$  since 
$\alpha_i \in (0,1)$ for $i=\overline{1,n-1}$

\begin{lemme}
Let condition {\bf A} be fulfilled, let all the functions   $K_i(t,s), i=\overline{1,n}$
be differentiable wrt $t$ and continuous wrt $s.$  Then the homogenius equation 
\begin{equation}
\int\limits_0^t K(t,s) s^N u(s) ds =0 \label{eq5}
\end{equation}
has  the trivial solution  in the space 
$\mathbb{C}_{[0,T]}.$ 
\end{lemme}

\proof
Let us differentiate the equation (\ref{eq5}) and take into account  (\ref{eq2}).
Then we  get an equivalent integral-functional equation
\begin{equation}
L u + \sum\limits_{i=1}^n \int\limits_{\alpha_{i-1}t}^{\alpha_i t} \frac{K_i^{\prime}(t,s)}{K_n(t,t)} \biggl(\frac{s}{t}\biggr)^N u(s) ds =0,
\label{eq6}
\end{equation}
where
$$ Lu = \sum\limits_{i=1}^n \biggl( \alpha_i^{1+N} K_i(t,\alpha_i t) u(\alpha_i t) -$$
$$- \alpha_{i-1}^{1+N} K_i(t,\alpha_{i-1} t) u(\alpha_{i-1} t) \biggr) (K_n(t,t))^{-1}. $$
Due to the condition {\bf A} in the space $\mathbb{C}_{[0,T]}$ we have the following 
estimate
$$ ||Lu-u  || \leq q ||u ||. $$
Therefore, according to the theorem on the inverse operator (\cite{13}, p.134), and
because of the inequalities
$0=\alpha_0 < \alpha_1 < \dots < \alpha_n =1$ there
exists the following bounded inverse operator
 $L^{-1} \in \mathcal{L}(\mathbb{C}_{[0,T]} \rightarrow \mathbb{C}_{[0,T]})$ 
 \begin{equation}
 || L^{-1} || \leq \frac{1}{1-q}
 \label{eq7}
 \end{equation}
 and equation (\ref{eq5}) can be reduced as follows 
 \begin{equation}
 u(t) = -L^{-1} \sum\limits_{i=1}^n \int\limits_{\alpha_{i-1}t}^{\alpha t} \frac{K_i ^{\prime}(t,s)}{K_n(t,t)} (s/t)^N u(s) ds \equiv A u,
 \label{eq8}
 \end{equation} where
 $0 \leq t \leq T. $
  Let us introduce the following equivalent norm  $||u|| = \max\limits_{0\leq t \leq T} e^{-lt} |u(t)|, \, l>0$  in the space $\mathbb{C}_{[0,T]}.$
  In this norm the inequality (\ref{eq7}) remains correct and for sufficiently
large $l$ an operator $A$ will be contracting since
  $||A|| \leq \frac{1}{1-q} m(l),$ where $m(l) \rightarrow 0$ for  $l\rightarrow +\infty.$
  Therefore homogeneous equation (\ref{eq8}) has  trivial solution.

  \endproof
   
  \begin{cor}
  Let all the condition of Lemma 1 be fulfilled,  $g(t) \in \mathbb{C}_{[0,T]}^{(1)}, \, |g^{\prime}(t)| ={o}(t^N)$ for $t\rightarrow +0.$
  Then inhomogeneous equation
  $\int\limits_0^t K(t,s) s^N u(s) ds = g(t)$
  has  the unique solution, and $u(t) \rightarrow 0$ for $t \rightarrow +0.$
  \end{cor}

  \textup{
Proof is trivial since differentiation of this equation leads to the 
equivalent equation
  \begin{equation}
  u(t) = Au + t^{-N} L^{-1} g^{\prime}(t)
  \label{eq9_1}
  \end{equation}
  with contracting operator $A$ and continuous  free function.
 }
  \begin{theorem}
  Let in the space  of continuous on $(0,T]$ functions which
have the finite limit for $t \rightarrow +0$ (briefly, in class  $\mathbb{C}_{(0,T]}$) 
 exists the function $x^N(t)$ such as for $t \rightarrow +0$
  $$\biggl (-\int\limits_0^t K(t,s) x^N(s) ds + f(t)   \biggr)^{\prime} = {o}(t^N),$$
$f(t) \in \mathbb{C}_{[0,T]}^{(1)},$ $f(0)=0.$
  Then equation (\ref{eq1})  has the following solution
  \begin{equation}
  x(t) = x^N(t) + t^N u(t)
  \label{eq10}
  \end{equation}
in class $\mathbb{C}_{(0,T]}.$
 Here function $u(t) \in \mathbb{C}_{[0,T]},$ $u(t) \rightarrow 0$ for $t \rightarrow +0$ and it
 can be uniquely constructed with successive approximations method.
 \end{theorem}

\proof
Proof follows from the corrolary 1. Indeed, with  (\ref{eq10})  we can rewrite the equation (\ref{eq1}) as follows
\begin{equation}
\int\limits_0^t K(t,s) s^N u(s) ds = g(t),
\label{eq11}
\end{equation}
 where function $g(t)$ is following
 \begin{equation}
 g(t) = -\int\limits_0^t K(t,s) x^N(s) ds + f(t) 
 \label{eq12}   
 \end{equation}
and it satisfies the condition of the corollary 1. Therefore in  (\ref{eq10})
the function $u(t)$ can be uniquely constructed with successive 
approximations from the equation (\ref{eq9_1}) using arbitrary initial condition.
\endproof
\begin{defi}
The equation (\ref{eq11}) with right hand side (\ref{eq12})  
we call {\it regularization} of the equation  (\ref{eq1}),
and function $x^N(t)$ as asymptotic approximation of solution (\ref{eq10}) of the equation (\ref{eq1}).
\end{defi}

\textup{
It is to be noted that one could numerically find the function $u(t)$ by solution of the 
equation (\ref{eq11}) based on well known numerical quadrature schemes
(see e.g. the bibliography in the monograph \cite{2}).
The method of constructing asymptotic approximations $x^N(t)$ 
in the solution (\ref{eq11}) we will study below in the Section 2.}

\section{The method of  asymptotic approximations construction }

\textup{
Let us suppose that along with the condition {\bf A} the  condition
{\bf B} be fulfilled. Functions $K_i(t,s), i=\overline{1,n}, \, f(t)$ are  $N+1$ times 
differentiable in the neighborhood of zero, where $N$ 
is selected according to the condition  {\bf A.}
We introduce an auxiliary algebraic equation wrt $j \in \mathbb{N}$}
\begin{equation}
L(j)  \triangleq  \sum\limits_{i=1}^n K_i(0,0) (\alpha_i^{1+j} - \alpha_{i-1}^{1+j}) = 0
\label{eq13}
\end{equation}
\textup{and name it as {\it characteristic} equation of  the integral equation  (\ref{eq1}).
Since $f(0)=0,$ then equation}
$$\sum\limits_{i=1}^n (\alpha_i K_i(t, \alpha_i t)x(\alpha_i t)) - \alpha_{i-1} K_i(t, \alpha_{i-1} t)x(\alpha_{i-1} t)))+ $$
$$+\sum\limits_{i=1}^n \int\limits_{\alpha_{i-1}t}^{\alpha_it} K_i^{\prime} (t,s) x(s) ds = f^{\prime}(t)$$
\textup{is equivalent to the equation (\ref{eq1}).
We will look for the asymptotical approximation of it's solution as following polynomial
$x^N(t) = \sum\limits_{j=0}^N x_j(\ln t) t^j.$}
\textup{Based on the method of undetermined coefficients, and taking into account the inequalities $0=\alpha_0 < \alpha_1 < \dots < \alpha_n =1$
we construct the recursive sequence of difference equations wrt  the coefficients 
$x_j(z) (z=\ln t)$ as follows:}
$$
K_n(0,0)x_j(z) + \sum\limits_{i=1}^{n-1} \alpha_i^{1+j} \bigl (K_i(0,0) - 
$$
\begin{equation}
-K_{i+1}(0,0)\bigr) x_j(z+a_i)= M_j(x_0,\dots, x_{j-1}),
\label{eq14}
\end{equation}
\textup{where $j=\overline{0,N}, \, a_i =\ln \alpha_i, \, i=\overline{1,n-1},$ $M_0=f^{\prime}(0).$\\
Here we follow (\cite{14}, p.330) and 
we seek the solution of the homogeneous difference equations in the form of
$x=\lambda^z.$}
\textup{Substitution of the function $\lambda^z$ into the homogenius 
difference equations leads to   $N+1$   equations for 
difference equations (\ref{eq14}):}
\begin{equation}
\mathcal{P}_j(\lambda) \equiv K_n(0,0) + \sum\limits_{i=1}^{n-1} \alpha_i^{1+j} \bigl(K_i(0,0) - 
\label{eq15}
\end{equation}
$$
K_{i+1}(0,0)\bigr) \lambda^{a_i}=0, \, j=\overline{0,N}.
$$
\textup{Therefore we have}\\
{\bf Property 1.} 
 $j$th equation (\ref{eq15}) has the root $\lambda=1$ if and only if $j$ satisfies
characteristic equation (\ref{eq13}) of the integral equation (\ref{eq1}). 
Moreover, multiplicity of the root $j$ of the equation (\ref{eq13}) is equal to $r_j$ 
iff
\begin{equation}
L(j) = \sum\limits_{i=1}^n K_i(0,0) (\alpha_i^{1+j} - \alpha_{i-1}^{1+j}) = 0,
\label{eq16}
\end{equation} 
\begin{equation}
\sum\limits_{i=1}^{n-1} \alpha_i^{1+j} (K_i(0,0)-K_{i+1}(0,0)) a_i^l =0, l=\overline{1,r_j-1},
\label{eq17}
\end{equation}
$$
\sum\limits_{i=1}^{n-1} \alpha_i^{1+j} (K_i(0,0)-K_{i+1}(0,0)) a_i^{r_j} \neq 0,
$$
where $\alpha_0 =0, \alpha_n=1, a_0 =0, a_n=0, a_i=\ln \alpha_i, i=\overline{1,n-1},$
and multiplicity    $r_j\leq n-1.$

\textup{Proof follows from the equality 
\begin{equation}
\sum\limits_{i=1}^{n-1} \alpha_i^{1+j} K_{i+1}(0,0) =  \sum\limits_{i=2}^n \alpha_{i-1}^{1+j} K_{i}(0,0)  
\label{eq18}
\end{equation}
and from the structure of the equations (\ref{eq13}), (\ref{eq15}).}
 
\textup{If we suppose that for certain $j$ multiplicity $r_j\geq n,$
then $K_1(0,0) = K_2(0,0) = \dots$ $\dots = K_{n-1} (0,0) = K_n(0,0)$
due to (\ref{eq17}),  since  $\det || a_i^l ||_{i,l=\overline{1,n}} \neq 0.$
But due to (\ref{eq16}) $K_n(0,0)=0,$
which contradicts {\bf A.} }

\textup{Under the conditions {\bf A, B} there are two cases.}

\subsection{Regular case} 

\textup{Let $L(j)\neq 0, j\in \mathbb{N}.$
Then $\lambda =1$ does not satisfy any of the  equations in the sequence (\ref{eq15}).
All the coefficients  $x_i$ of the asymptotics $x^N = \sum\limits_{i=0}^N x_i t^i $
can be determined uniquely with method of undetermined coefficients and do not depend
upon $\ln t$.}

\textup{Therefore we have the following theorem }
\begin{theorem}
     Let the conditions {\bf A}, {\bf B} and $L(j) \neq 0, j\in \mathbb{N}$
be fulfilled. Then equation (\ref{eq1})
   has in $\mathbb{C}_{[0,T]}$ the solution  $x(t) = \sum\limits_{i=0}^N x_i t^i + t^N u(t),$
   where  $x_i $ are determined uniquely with method of undetermined coefficients and function
    $u(t)$ is uniquely constructed (numerically or with successive approximations
from equation (\ref{eq11})).
\end{theorem}

 \subsection{Irregular case}
 
 \textup{Let $L(j)=0$ only for $j \in \{  j_1, \dots , j_k \} \subset \{ 0, 1, \dots  , N \} $
  and multiplicity of the root $\lambda =1$  for the corresponding 
characteristic equation  is  $r_j.$
  Let in $j$th difference equation (\ref{eq14}) right hand side $M_j(z)$ 
appears to be polynomial from  $z$
  of the order $n_j \geq 0.$
   Then in irregular case, i.e. for $r_j \geq 1,$
 based on (\cite{14}, p.338)
particular solution of  the $j$th equation (\ref{eq14}) we have to search as following
polynomial
 $\hat{x}(z) = \sum\limits_{i=rj}^{n_j + r_j}  c_i z^i.$
 Coefficients $c_i$ of this polynomial can be sequentially calculated by the method of undetermined coefficients  starting from $c_{n_j+r_j}.$  Coefficient  $x_j(z)$  of desired asymptotical
approximation
 $x^N$ in this case is as follows $$x_j(z) =  c_0 + c_1 z + \dots + c_{r_{j-1}} + \hat{x}(z). $$
 In irregular case when $r_j \geq 1,$
constants $c_0, \dots, c_{r_j-1}$ remain arbitrary since functions $z^i, \, i=0,1, \dots, r_j-1$
 satisfy  $j$th homogenius difference equation corresponding to (\ref{eq14}).}
 
\textup{ In applications one could use the Propetry 1, and
coefficient $x_j(z)$ in irregular case directly as polynomial
 $\sum\limits_{i=0}^{n_j + r_j} c_i z^i,$
where $c_{n_j + r_j}, \dots , c_0$ are determined sequentially 
using method of undetermined coefficients. And
  $c_{r_j-1}, \dots , c_0$ remains arbitrary. Therefore
  in irregular case when  $L(j)=0$ for some $j$ new arbitrary constants $r_j$
appear in determination of the coefficient $x_j(z).$ Order of the polynomial $x_j(z)$ on
the value of the multiplicity of $r_j$ or root
     $\lambda=1$  of the $j$th equation (\ref{eq15}) becomes greater than order
     $n_j$ of the right hand side of the corresponding equation  (\ref{eq14}), i.e. of the order of the polynomial $M_j(z).$}

    \textup{Therefore we have the following theorem:}
    \begin{theorem}
    Let the conditions {\bf A, B} are fulfilled. Let characteristic equation $L(j)=0$ of integral
equation (\ref{eq1}) has exactly 
     $k$ natural roots $\{ j_1,\dots , j_k \}.$ And let the root $\lambda=1$ of  $j$th
     equation (\ref{eq15})  has multiplicity $r_j.$
    Then equation (\ref{eq1}) has the following solution in $\mathbb{C}_{(0,T]}$
\begin{equation}
x =\sum\limits_{i=0}^N x_i(\ln t) t^i + t^N u(t),
\label{eq18a}
\end{equation}
    which depends on $p=r_1 + \dots + r_k$ arbitrary constants. Moreover,
coefficints  $x_i$ of the asymptotic approximation
$x^N(t)$ are polynomials from $\ln t.$ 
     Function $u(t)$ can be
      constructed by successive approximations which converge uniformly for  $t \in [0,T],$
    or numerically from (\ref{eq11}).    
    \end{theorem}
     \begin{rem}  If $L(0)=0,$ then in solution (\ref{eq18a}) $x_0 = const + a \ln t,$
    where $a$ is the defined constant. Therefore in this case $x(t) \in \mathbb{C}_{(0,T]}, \, \lim\limits_{t\rightarrow +0} x(t) = \infty.$
    \end{rem}
     \begin{rem} These results can be generalized if in the equation (\ref{eq1})
 for $\alpha_{i-1}(t)  \leq s \leq \alpha_i(t)$
$K(t,s) = K_i(t,s)$ 
where  $\alpha_i(0) =0, $  $0 \leq \alpha_0^{\prime}(0) < \alpha_1^{\prime}(0) < \cdots < \alpha_n^{\prime}(0) \leq 1,$
$i=1,\dots , n.$
     \end{rem}

\section*{Conclusion}
\textup{Our method 
allow studies of the Volterra integral equations of the first kind
when kernals are discontinuous operator-functions
operating in Banach spaces. Our work naturally complements the theory of  
integral models of evolving systems.}

\section*{Acknowledgement}
\textup{
This work is fulfilled as part of Federal Framework Programm
 ``Cadri''  P696 (30.09.2010), it is also
partly supported by RFBR, projects \mbox{ 09--01--00377}
and \mbox{ 11-08-00109}.}

\bibliographystyle{alpha}
\bibliography{./exemple-biblio}

\begin{thebibliography}{MST11}

\bibitem[Apa03]{2}
A.S. Apartsyn.
\newblock {\em Nonclassical Linear Volterra Equations of the First Kind}.
\newblock (first edition in 2003 by the Walter de Gruyter ), 2003.

\bibitem[DL95]{3}
A.M. Denisov and A.~Lorenzi.
\newblock On a special volterra integral equation of the first kind.
\newblock {\em Boll. Un. Mat. Ital. B.}, 7(9):443--457, 1995.

\bibitem[Eva10]{6}
G.C. Evans.
\newblock Integral equation of the second kind with discontinuous kernel.
\newblock {\em Transactions of the American Mathematical Society},
  11(4):393--413, 1910.

\bibitem[Gel59]{14}
A.O. Gelfond.
\newblock {\em Calculus of Finite Differences}.
\newblock Fizmatlit, 1959.

\bibitem[HY96]{1}
N.~Hritonenko and Yu. Yanenko.
\newblock {\em Modeling and Optimization of the Lifetime of Technologies}.
\newblock Kluwer Academic Publushers, 1996.

\bibitem[HY03]{5}
N.~Hritonenko and Yu. Yatsenko.
\newblock {\em Applied mathematical modelling of engineering problems}.
\newblock Kluwer, 2003.

\bibitem[Khr06]{7}
A.P. Khromov.
\newblock Integral operators with discontinuous kernel on piecewise linear
  curves.
\newblock {\em Sbornik: Mathematics}, 197(11):115--142, 2006.

\bibitem[Mag83]{8}
N.A. Magnitsky.
\newblock Asymptotics of the solution of the volterra integral equations of the
  first kind.
\newblock {\em DAN USSR}, 169(1):29--32, 1983.

\bibitem[MST11]{4}
E.V. Markova, I.V. Sidler, and V.V. Trufanov.
\newblock On models of developing systems and their applications.
\newblock {\em Automation and Remote Control}, 72(7):1371--1379, 2011.

\bibitem[SS06]{9}
N.A. Sidorov and D.N. Sidorov.
\newblock Generalized solutions to integral equations in the problem of
  identification of nonlinear dynamic models.
\newblock {\em Differential Equations}, 42(9):1312--1316, 2006.

\bibitem[SS11]{12}
N.A. Sidorov and D.N. Sidorov.
\newblock Small solutions of nonlinear differential equations near branching
  points.
\newblock {\em Russian Mathematics}, 55(5):43--50, 2011.

\bibitem[SST07]{11}
N.A. Sidorov, D.N. Sidorov, and A.V. Trufanov.
\newblock Existence and structure of solution of integral-functional volterra
  equations of the first kind.
\newblock {\em Bulletin of Irkutsk State University: Mathematics},
  1(1):267--274, 2007.

\bibitem[ST09]{10}
N.A. Sidorov and A.V. Trufanov.
\newblock Nonlinear operator equations with a functional perturbation of the
  argument of neutral type.
\newblock {\em Differential Equations}, 45(12):1840--1844, 2009.

\bibitem[Tre07]{13}
V.A. Trenogin.
\newblock {\em Functional Analysis}.
\newblock Nauka, Moscow, 2007.

\end{thebibliography}


\end{document}